\documentclass[]{interact}

\usepackage{epstopdf}
\usepackage[caption=false]{subfig}

\theoremstyle{plain}

\theoremstyle{definition}

\theoremstyle{remark}

\usepackage{graphicx}
\usepackage{amsmath,amsthm,amssymb,amsopn,amsfonts}
\usepackage{mathtools}
\usepackage{multirow}
\usepackage{bm}

\usepackage[ruled,vlined]{algorithm2e}

\renewcommand{\vec}[1]{\textbf{#1}}

\usepackage{multirow}

\usepackage{todonotes}

\usepackage{capt-of}

\usepackage{hyperref}

\begin{document}

\title{Structural Gaussian Priors for Bayesian CT reconstruction of Subsea Pipes}

\author{
\name{Silja L. Christensen\thanks{CONTACT Silja L. Christensen. Email: swech@dtu.dk}, Nicolai A. B. Riis, Felipe Uribe, and Jakob S. Jørgensen}
\affil{Department of Applied Mathematics and Computer Science, Technical University of Denmark, Kongens Lyngby, Denmark}
}

\maketitle

\begin{abstract}
A non-destructive testing (NDT) application of X-ray computed tomography (CT) is inspection of subsea pipes in operation via 2D cross-sectional scans. Data acquisition is time-consuming and costly due to the challenging subsea environment. While reducing the number of projections in a CT scan yields time and cost savings, this compromises the reconstruction quality when conventional methods are used. To address this issue we take a Bayesian approach to CT reconstruction and focus on designing an effective prior that introduces additional information to the problem. We propose a new class of \emph{structural Gaussian priors} to enforce expected material properties in different regions of the reconstructed image based on available a-priori information about the pipe's layered structure. The prior is composed from Gaussian distributions such that we can define an efficient computational strategy to sample the resulting posterior distribution. This is essential in practical NDT applications for fast processing of the large-scale data we see in CT. Numerical experiments with synthetic and real data show that the proposed structural Gaussian prior reduces artifacts and enhances contrast in the reconstruction, compared to using only a conventional Gaussian Markov Random Field (GMRF) prior or no prior information at all.
\end{abstract}

\begin{keywords}
X-ray computed tomography; Uncertainty Quantification; Bayesian Inversion; Non-destructive testing; Pipeline inspection
\end{keywords}

\section{Introduction}

Subsea pipelines are used to transport oil and gas around the world. It is crucial that the subsea pipes are in good condition to avoid environmentally damaging leaks. Therefore, efficient methods for non-destructive inspection are in demand. One suitable method is X-ray computed tomography (CT) that is used to obtain 2D cross-sectional images of the pipes using data from a finite number of projections \cite{Hansen2021}. To produce high-quality images with limited noise and artifacts, standard CT reconstruction methods require lengthy measurements at many projection (view) angles. To reduce costs, it is desirable to limit measurement time, and this can be achieved by reducing the number of view angles. However, if standard iterative reconstruction methods are used, this compromises the reconstruction quality significantly \cite{hansen_2010}. Moreover, since the data is typically scarce and noisy, there is substantial uncertainty affecting the reconstructions. Therefore, we follow a Bayesian approach, where the idea is to update a \emph{prior} probability distribution of the image by including information about the CT model and projection data to obtain a so-called \emph{posterior} distribution describing the reconstruction and its uncertainty. Our idea is to use prior knowledge of the pipes' internal structure to improve the image reconstruction quality, even for cases with measurements taken at only few view angles.

Previous work has taken a deterministic approach to CT scanning of subsea pipes \cite{Riis2018}. The authors use microlocal analysis to design a favorable offset scan geometry and propose a reconstruction method based on shearlets with a weighted sparsity penalty. However, this method does not account for and quantify the uncertainty in the reconstructions, and thus a Bayesian approach is needed for a reliable solution of the CT problem. 

In Bayesian CT reconstruction \cite{Kaipio2005, Suuronen2020, Uribe2021}, it is important to choose an appropriate prior. Typical choices in this setting are Markov Random Field (MRF) priors based on multivariate Gaussian and Laplacian distributions. These represent a statistical parallel to the well-known Tikhonov and total variation (TV) regularization methods, which are commonly used in inverse problems and in particular CT \cite{Hansen2021,sidky2008image}. An overview of several MRF type priors applied to 2D inverse problems (including CT) is given in \cite{Bardsley2018} and more detailed descriptions can be found in earlier papers \cite{Bardsley2013, Bardsley2012a}. A recent approach is the Cauchy differences prior (also a MRF-type prior), which is edge-preserving and discretization invariant \cite{Markkanen2019,suuronen_et_al_2022}. All MRF priors have regularizing effects, and depending on the scanned object, a smoothness-promoting or edge-preserving prior can be chosen. However, a common characteristic for these distributions is that they do not model specific features, such as the layer structure in subsea pipes, when used as priors in Bayesian reconstruction. Furthermore, non-Gaussian MRF-type priors require more advanced sampling techniques, which in turn are also computationally expensive and can therefore be prohibitive for large-scale problems.

If a particular structure of the target image is known in advance, a stronger prior can be defined. A series of studies use this idea by forming a texture-preserving prior for low-dose CT reconstruction using a previous normal-dose scan \cite{Zhang2014,Zhang2016,Han2016,Gao2019}. Another study exploits the fact that the human body consists of a finite number of tissue types with well-known attenuation coefficients and the study designs the prior accordingly including an inherent classification problem \cite{Fukuda2010}.

The driving idea in this work is to utilize as much prior information as possible. We know the overall layout of the pipes because they are manufactured according to specific design standards. This leads to our proposed \textit{structural Gaussian prior} that introduces information about the pipes' layered structure as well as their constitutive materials. Numerical experiments show that this prior improves the reconstruction quality compared to standard priors, especially if only few view angles are available. Moreover, the Bayesian approach allows probabilistic modelling of noise and prior information, and through sampling of the resulting posterior distribution, we can obtain reconstructions with uncertainty estimates. Another advantage of the structural Gaussian prior is that it can be sampled efficiently in large-scale problems. 

Our paper is organized as follows. In Section \ref{sec:BayesianCT} we introduce the Bayesian approach to CT reconstruction. Section \ref{sec:SGP} describes our proposed structural Gaussian prior. In Section \ref{sec:post_sample} we derive the posterior and how to sample it efficiently. In Section \ref{sec:experiments} we perform numerical experiments with synthetic and real data, and in Section \ref{sec:discussion} we discuss and conclude our findings.

\section{Bayesian CT problem} \label{sec:BayesianCT}

\begin{figure}
    \centering
    \raisebox{-0.5\height}{\includegraphics[width = 0.3\textwidth,angle=-90,origin=c]{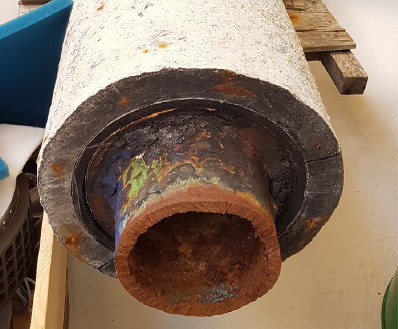}}
    \hfill
    \raisebox{-0.5\height}{\includegraphics[width = 0.68\textwidth]{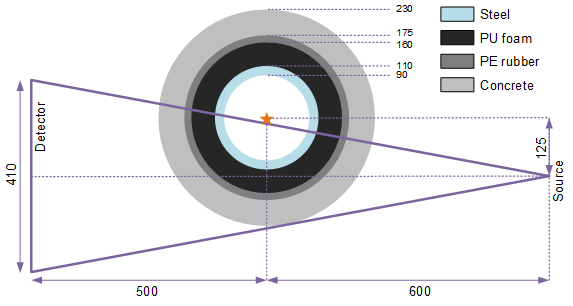}} 
    \caption{\textit{Left:} Photograph of subsea oil pipe. \textit{Right:} Pipe and offset scan geometry specifications; dimensions in mm. The CT rotation axis is marked with a star.}
    \label{fig:scan_geom} 
\end{figure}

Deep sea oil pipes are large and consist of dense materials including steel and concrete. In order to penetrate the pipe with X-rays, it is necessary to use high power and a narrow beam source. The narrow beam cannot illuminate the full pipe, but an exterior tomography setup can be used. In this case, the X-ray source is offset from the pipe centre, which is also the rotation axis. It can be shown that this setup is ideal, if we are interested in recovering the pipe layer structure and defects therein rather than the center of the pipe \cite{Riis2018}. We use real data acquired from a prototype set-up of an oil pipe CT scanner. A linear detector panel with 510 sensors of sizes 0.8 mm each was used, and 360 equi-angular projections were recorded in a full 360 degree rotation. The pipe, its specifications, and the offset fan-beam acquisition geometry are illustrated in Figure \ref{fig:scan_geom}.

\subsection{CT model}

If we discretize the reconstruction domain in pixels, the inverse problem associated with CT can be expressed as a linear system of equations,
\begin{equation}
    \vec{d} = \vec{Ax} + \vec{e}, \qquad \vec{e}\sim\mathcal{N}\left(\vec{0}, \frac{1}{\lambda}\vec{I}_n\right), \label{eq:forward}
\end{equation}
where $\vec{d} \in \mathbb{R}^m$ is the observed X-ray absorption (known as a sinogram), $\vec{x}\in \mathbb{R}^n$ is the unknown vector of linear attenuation coefficients in the reconstruction domain with $N\times N = n$ pixels, $N$ is the domain discretization size, and $\vec{A}\in \mathbb{R}^{m \times n}$ contains entries $a_{ij}$ that describe the intersection of ray $i$ with pixel $j$ (see, e.g., \cite{Hansen2021} for details). We model the data noise as additive Gaussian, suitable for the high-intensity X-rays used in the data collection process. That is, $\vec{e}\in \mathbb{R}^m$ is independent and identically distributed Gaussian noise with mean zero and precision parameter $\lambda$, and $\vec{I}_n \in \mathbb{R}^{n\times n}$ is an identity matrix. 
 
\subsection{Bayesian inversion}

We take a Bayesian approach to the CT problem. Using Bayes' Theorem we can characterize the solution in terms of the posterior distribution given by the probability density function (PDF),
\begin{equation}
    \pi_\text{post}(\vec{x}\,|\,\vec{d}) \propto \pi_\text{pr}(\vec{x})\pi_\text{li}(\vec{d}\,|\,\vec{x}). \label{eq:bayes}
\end{equation}

Here, $\pi_\text{pr}(\vec{x})$ denotes the prior PDF describing some known properties of the attenuation coefficients $\vec{x}$. The prior is modelled as a Gaussian distribution with PDF given by
\begin{equation}
	\pi_\text{pr}(\vec{x}) = \frac{1}{\left(2\pi\right)^{n/2}}\left(\left|\vec{R}^T\vec{R}\right|\right)^{1/2} \exp\left(-\frac{1}{2}\lVert \vec{R}(\vec{x}-\bm{\mu}) \rVert_2^2\right), \label{eq:Gaussian}
\end{equation}
where $\bm{\mu} \in \mathbb{R}^n$ is the mean and $(\vec{R}^T\vec{R})\in \mathbb{R}^{n\times n}$ the precision matrix with full rank. We refer to $\vec{R}$ as the square-root precision matrix, and $|\cdot|$ denotes the matrix determinant. Below, in Section \ref{sec:SGP}, the specifics of how we build the pipe structure into the Gaussian prior are given.

From \eqref{eq:forward}, we obtain the distribution of the data $\vec{d}$ given some known vector of linear attenuation coefficients $\vec{x}$,
\begin{equation}
	\pi_\text{li}(\vec{d}\,|\,\vec{x}) = \left(\frac{\lambda}{2\pi}\right)^{m/2} \exp\left( - \frac{\lambda}{2} \left\Vert \vec{A}\vec{x}-\vec{d} \right\Vert_2^2\right), \label{eq:likelihood}
\end{equation}
which can be seen as a likelihood function of $\vec{x}$ given some observed data $\vec{d}$.

Because the prior and likelihood are Gaussian, we can obtain a closed-form expression for the posterior, which is also Gaussian (see Section \ref{sec:postdistb}). In general, it is numerically impractical to compute the posterior precision matrix, although the posterior mean can be simple to obtain. This is because the distribution is Gaussian, and thus the mean is equivalent to the maximum a-posteriori probability (MAP) estimate, which is often used as a point estimator to describe the posterior, as it can be easily computed through optimization. Nevertheless, the MAP estimator itself does not provide any uncertainty information, and therefore does not reach the full potential of the Bayesian approach. Indeed, one of the strong points of Bayesian inversion is the ability to give uncertainty estimates. Therefore, our approach entails a method for sampling the posterior distribution as described in Section \ref{sec:sampling}, instead of focusing on a single point estimate.

\section{Structural Gaussian priors}\label{sec:SGP}

We must specify a prior in order to perform Bayesian CT reconstruction. We work with Gaussian priors because they have nice mathematical properties making them easy to combine, and because they can be sampled efficiently compared to more general distribution models. First, we outline the prior structural information that we have. Then we describe our proposed method of incorporating the known pipe geometry into what we call \textit{structural Gaussian priors} (SGPs). The SGP combines (stacks) priors for each of the materials in the pipe, thereby enforcing a different prior in each region. In addition to the pipe material and geometry information, the SGP also contains a prior on the pixel differences of the whole image. The purpose of this is to add a global smoothness assumption on the image, which also ensures full rank of the precision matrix of the combined Gaussian prior in \eqref{eq:Gaussian}.

\subsection{Prior structural information}\label{sec:structurepriors}

We have prior information about the pipe's layered structure and constitutive materials, which we will incorporate in the prior probability model. This is the core idea of our proposed SGP. To build the SGP we need: 1) masks for each pipe layer and for the background, and 2) estimates of the linear attenuation coefficients of the layer materials. These are described in the following subsections.

\begin{figure}
    \centering
    \includegraphics[width=0.4\textwidth, trim={2.5cm 0.9cm 0cm 1.1cm},clip]{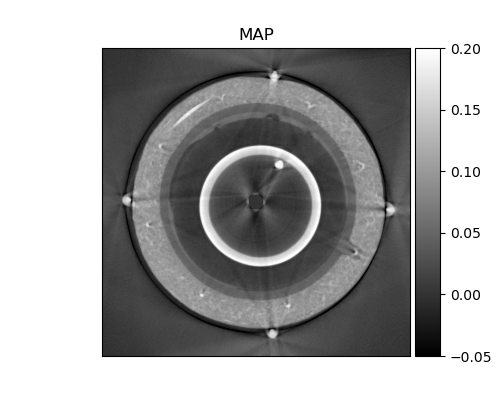}
    \includegraphics[width = 0.55\textwidth]{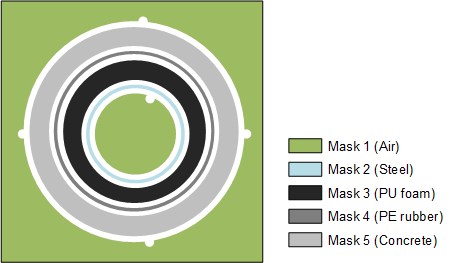}
    \caption{\textit{Left:} Real pipe with four layers of different materials and a background of air. \textit{Right:} SGP masks for $i=1,\ldots,5$ that each represent a material with attenuation coefficient $\alpha_i$ given in \eqref{eq:abscoeff}. The white regions do not have any prior promoting specific attenuation coefficients. The index $i=0$ is used to denote the prior on the whole image domain (background and pipe).}
    \label{fig:SGPmask_real}
\end{figure}

\subsubsection{Pipe layer masks}

Given $p$ non-overlapping regions of different materials, we represent each region by a mask, i.e., a binary vector $\vec{m}_i\in \mathbb{R}^n$ for $i=1,\ldots,p$, in which the value of element $j$ is 1, if pixel $j$ is in region $i$ and $0$ if not, for $j=1,\ldots,n$. Figure \ref{fig:SGPmask_real} shows a reconstruction of the real pipe and the corresponding masks (with $p=5$). We use masks that are slightly smaller than the actual pipe layers to account for modelling error in the pipe structure and scan geometry. The masks are constructed based on a high-quality standard CT reconstruction of the real pipe. This may also be obtained from lower-quality data or from available CAD drawings. We further define mask matrices $\vec{M}_{i} \in \mathbb{R}^{l_i\times n}$, $i=1,\ldots,p$, that pick out the $l_i$ pixels of region $i$ of the $n$-dimensional image. Specifically, $\vec{M}_i$ is formed by keeping only the $l_i$ rows of the identity matrix $\vec{I}_n$ corresponding to where there are 1-values in the mask $\vec{m}_i$.

\subsubsection{Linear attenuation coefficients of pipe materials} \label{sec:abscoeff}

When the X-ray energy is known, we can estimate the linear attenuation coefficients $\alpha_{i}$ (ignoring scattering) for the materials as \cite{buzug2008photon}:
\begin{equation}
    \alpha_i = 1/\ell_i = \kappa_i(E) \rho_i, \qquad i=1,\ldots,p, \label{eq:abscoeff}
\end{equation}
where $\rho_i$ is the material density and $\kappa_i(E)$ is the mass attenuation coefficient at X-ray energy $E$ that can be found, e.g., in the NIST database \cite{NIST}. Equation \eqref{eq:abscoeff} also indicates that the linear attenuation coefficient is sometimes given in terms of the mean free path of the material $\ell_i$.

For thick high-density materials, the scattering of X-rays cannot be ignored and it becomes necessary to account for a physical effect referred to as buildup in the material. Consider the X-ray intensity at some point in the material. Some of the intensity will be due to photons that have reached the point undisturbed (primary) and some of the intensity will be due to scattered photons (secondary). The buildup is defined at a given point as the ratio of the total (primary and secondary) intensity to the primary intensity. With buildup, both primary and secondary intensity reach the detector, which then gives a higher intensity reading compared to the case where buildup is not present, and only primary intensity reaches the detector. With a so-called buildup factor we can account for the effect. 

Assuming a homogeneous material, the X-ray attenuation including buildup is described as \cite{Rafiei2018}:
\begin{equation}
    I = I_0 \, B(E, w_i/\ell_i) \exp\left(-w_i/\ell_i\right),  \label{eq:lamberbeer_extended}
\end{equation}
where $I_0$ is the X-ray intensity entering the homogeneous material, $I$ is the X-ray intensity leaving the homogeneous material, $w_i$ is the width of the homogeneous material $i$, and $B(E, w_i/\ell_i)$ is the buildup factor which can be found from e.g. \cite{buildup}. However, in most CT implementations, the forward operator $\vec{A}$ is constructed ignoring buildup and assumes that the X-ray attenuation follows Lambert--Beer's law:
\begin{equation}
    I = I_0 \exp\left(-\alpha_i w_i\right). \label{eq:lambertbeer}
\end{equation}

Using Lambert--Beer's law for CT reconstruction of dense materials is essentially a model error, and it will lead to linear attenuation coefficients that are lower than expected according to \eqref{eq:abscoeff} because the secondary intensity is measured but not included in the model. Instead of correcting the model error, we adjust our prior expectations about the linear attenuation coefficient to include the buildup factor. By equating the right-hand sides of \eqref{eq:lamberbeer_extended} and \eqref{eq:lambertbeer}, and inserting \eqref{eq:abscoeff} instead of $1/\ell_i$, we estimate the linear attenuation coefficient including the buildup factor for a dense material as:
\begin{equation}
    \alpha_i = \kappa_i(E)\rho_i - \frac{\log\left(B(E, w_i\kappa_i(E)\rho_i)\right)}{w_i}. \label{eq:abscoeffbuildup}
\end{equation}

In table \ref{tab:abscoeff} we will use \eqref{eq:abscoeff} for prior information about concrete, polyurethane, and polyethylene, and we will use \eqref{eq:abscoeffbuildup} for prior information about the much denser material, steel. 

\subsection{Structural Gaussian priors}\label{sec:StandardGaussians}

To form the SGP, we use Gaussian distributions defined in terms of means $\bm{\mu}_i\in\mathbb{R}^{n}$ and a symmetric positive semi-definite precision matrices $\left(\vec{R}_i^T\vec{R}_i\right)\in\mathbb{R}^{n\times n}$. In \eqref{eq:Gaussian} we gave an expression for Gaussians with precision matrices of full rank. When forming the SGP, we will need intrinsic Gaussian distributions that can be defined for rank-deficient precision matrices of rank $k$, where $n\geq k>0$ . The intrinsic Gaussian prior density is defined as \cite[Ch. 3]{RueHeld2005}:
\begin{equation}
	\pi_i(\vec{x}) = \frac{1}{\left(2\pi\right)^{k/2}}\left(\left|\vec{R}_i^T\vec{R}_i\right|^{*}\right)^{1/2} \exp\left(-\frac{1}{2}\lVert \vec{R}_i(\vec{x}-\bm{\mu}_i) \rVert_2^2\right), \label{eq:iGMRF}
\end{equation}
where $|\cdot|^{*}$ is a generalized determinant given by the product of the $k$ non-zero eigenvalues of the precision matrix. Note this density is multivariate Gaussian if and only if $k=n$. In the SGP, we use two different Gaussian priors that are specified below: 
\begin{description}
    \item[IID:] The independent and identically distributed (IID) Gaussian prior is a noise-dampening prior that favors posterior pixel values around the prior mean. We use it to express prior information of expected linear attenuation coefficients $\alpha_i$ in the $p$ different materials. The prior is defined as an intrinsic Gaussian parameterized by mean and square-root precision:
    \begin{align}
        \bm{\mu}_i = \alpha_i \vec{1} \qquad \text{and} \qquad \vec{R}_i = \sqrt{\delta_i} \vec{M}_i,
    \end{align}
    where $\vec{1}\in \mathbb{R}^{n}$ is a vector of ones and $\delta_i$ is a prior precision parameter. 
    \item[GMRF:] The Gaussian Markov Random Field (GMRF) prior is a smoothing prior that promotes neighbouring pixels to take the same value. The pipe materials can be considered roughly homogeneous and therefore the GMRF prior is enforced on the entire image domain. We index this prior with $i=0$ corresponding to the whole image, and in this case the mean and square-root precision are:
    \begin{align}
        \bm{\mu}_0 = \vec{0}, \qquad \text{and} \qquad \vec{R}_0 = \sqrt{\delta_0} \begin{bmatrix} \vec{D}_1\\ \vec{D}_2 \end{bmatrix},
    \end{align}
    where $\vec{0}\in \mathbb{R}^{n}$ is a vector of zeros, $\vec{D}_1 = \vec{I}_N \otimes \vec{D}$, $\vec{D}_2 = \vec{D} \otimes \vec{I}_N$, $\vec{D}\in\mathbb{R}^{(N+1) \times N}$ is a finite difference matrix with backward differences and zero Dirichlet boundary conditions, and $\otimes$ denotes the Kronecker product. Note that $\vec{R}_0$ has full column rank. Therefore the precision matrix has full rank, and the density is thus multivariate Gaussian.
\end{description}

To construct the SGP prior, we enforce IID priors on the masked image regions and a GMRF prior on the whole image domain as described above. Under the assumption that the priors $\pi_{i}(\vec{x})$, $i=0,\ldots,p$ are independent of each other, the SGP prior is given as a product of the individual priors:
\begin{equation}
    \pi_\text{pr}(\vec{x}) = \prod_{i=0}^p \pi_{i}(\vec{x}) \label{eq:prior}.
\end{equation}

Then, the SGP prior is also Gaussian with mean and square-root precision: 
\begin{equation}
\vec{R}_\text{pr} = \begin{bmatrix} \vec{R}_0 \\\vdots\\ \vec{R}_p \end{bmatrix} \qquad \text{and} \qquad \bm{\mu}_\text{pr} = \left(\vec{R}_\text{pr}^T\vec{R}_\text{pr}\right)^{-1}\sum_{i=0}^p \vec{R}_i^T\vec{R}_i \bm{\mu}_i. \label{eq:prior_meanprec}
\end{equation}

Note, that the full column rank of $\vec{R}_0$ leads to $\vec{R}_\text{pr}$ also having full column rank, and therefore the prior precision matrix $\vec{R}_\text{pr}^T\vec{R}_\text{pr}$ has full rank and is invertible. 

\subsubsection{SGP configurations}\label{sec:SGPconfig}

We use the pipe shown in Figure \ref{fig:SGPmask_real} as a case study for numerical experiments in Section \ref{sec:experiments}. The pipe consists of four different materials, steel in region 2, polyurethane (PU) foam in region 3, polyethylene (PE) rubber in region 4, and concrete in region 5. The pipe is scanned with air as the background and this will be denoted region 1. Thus we have $p=5$ masked regions where we can impose IID priors promoting the expected attenuation coefficients. We study two configurations of the SGP:
\begin{description}
    \item[SGP-BG:] Assuming we have information about the pipe's dimensions but not internal structure, we use a GMRF on the whole domain and add an IID prior on the air background (`BG') such that $\pi_\text{pr}(\vec{x})=\prod_{i=0}^1\pi_{i}(\vec{x})$.
    \item[SGP-F:] Assuming we also have information about the internal pipe structure, we again use a GMRF on the full (`F') domain and add IID priors on the background and on the pipe layers such that $\pi_\text{pr}(\vec{x})=\prod_{i=0}^5\pi_{i}(\vec{x})$.
\end{description}

\section{Gaussian posterior and efficient sampling} \label{sec:post_sample}

In Sections 2 and 3, we assumed Gaussian distributions for the likelihood and priors respectively, which in turn yields a Gaussian posterior PDF \cite{Kaipio2005}. This is mathematically convenient because we can derive a closed-form expression for this distribution and sample it efficiently. Below, we define the posterior PDF and then show how to sample it by solving a system of linear equations. 

\subsection{Posterior PDF} \label{sec:postdistb}

Using the likelihood in \eqref{eq:likelihood} and the prior in \eqref{eq:prior}-\eqref{eq:prior_meanprec}, the posterior PDF is given by Bayes' Theorem \eqref{eq:bayes} as:
\begin{align}
\pi_{\text{post}}(\vec{x}\,|\,\vec{d}) \propto{}& \left(\frac{\lambda}{2\pi}\right)^\frac{m}{2}\exp\left[-\frac{\lambda}{2} \left\Vert \vec{A}\vec{x}-\vec{d}\right\Vert_2^2 \right] \times \prod_{i=0}^p  \left(\left|\frac{\vec{R}_i^T\vec{R}_i}{2\pi}\right|^*\right)^\frac{1}{2} \exp\left[-\frac{1}{2}\lVert\vec{R}_i(\vec{x}-\bm{\mu}_i)\rVert_2^2\right]\nonumber\\
\propto{}& \exp\left[-\frac{\lambda}{2} \left\Vert \vec{A}\vec{x}-\vec{d}\right\Vert_2^2 - \frac{1}{2} \sum_{i=0}^p \lVert \vec{R}_i(\vec{x}-\bm{\mu}_i) \rVert_2^2\right]. \label{eq:GaussianPosterior0}
\end{align}

The above expression defines the Gaussian random vector
\begin{equation} \label{eq:GaussianPosteriorDistribution}
    \vec{x}\,|\,\vec{d} \sim \mathcal{N}(\bm{\mu}_\mathrm{post}, (\vec{R}^T_\mathrm{post}\vec{R}_\mathrm{post})^{-1}),
\end{equation}
where the posterior square-root precision is given by 
\begin{equation}
\vec{R}_\text{post} = \begin{bmatrix} \sqrt{\lambda}\vec{A}\\ \vec{R}_0 \\\vdots\\ \vec{R}_p \end{bmatrix} \in \mathbb{R}^{q \times n}  \qquad \text{with} \qquad q = \left(m+2N\left(N+1\right)+\sum_{i=1}^p l_i\right),\label{eq:postprec}
\end{equation}
and the posterior mean is given implicitly by the expression:
\begin{align}
    \vec{R}_\text{post}^T\vec{R}_\text{post}\bm{\mu}_\text{post} = \lambda \vec{A}^T\vec{d} + \sum_{i=0}^p \vec{R}_i^T\vec{R}_i \bm{\mu}_i.\label{eq:postmean0}
\end{align}

Since $\vec{R}_0$ and $\vec{R}_\text{post}$ have full column rank, the posterior precision is invertible, and the posterior distribution is Gaussian with well-defined covariance. That is, the closed-form expression for the posterior distribution is \eqref{eq:GaussianPosteriorDistribution}, which can be seen by factoring out \eqref{eq:GaussianPosterior0} and reordering the equation into the form \eqref{eq:Gaussian}.

The large-scale nature of CT makes it infeasible to form the matrix $\vec{A}$ explicitly, which can be avoided by matrix-free implementations of the application of $\vec{A}$ and its adjoint. Then, the posterior mean can be obtained by solving \eqref{eq:postmean0}. However, to compute the posterior precision, we must form $\vec{A}^T\vec{A}$, which is numerically impractical. Instead we sample the posterior as described below, and statistics of the samples are used to obtain an uncertainty estimate. 

\subsection{Posterior Sampling}\label{sec:sampling}

Our goal is to efficiently sample the posterior \eqref{eq:GaussianPosteriorDistribution} without explicitly forming the matrix $\vec{A}$. Inspired by \cite{Uribe2021}, we generate posterior realizations by solving the linear system of equations:
\begin{equation}
    \vec{R}_\text{post} \vec{x}^* = \begin{bmatrix} \sqrt{\lambda}\vec{d} \\ \vec{R}_0 \bm{\mu}_0 \\ \vdots\\ \vec{R}_p\bm{\mu}_p \end{bmatrix} + \bm{\xi}, \quad \bm{\xi}\sim \mathcal{N}(\vec{0},\vec{I}_{q}). \label{eq:post_linear_system}
\end{equation}
Note that the solution method must support matrix-free implementation of $\vec{A}$.

We can show that $\vec{x}^*$ is distributed according to the posterior by first forming the normal equations of the linear system above, i.e.
\begin{align}
\left(\lambda \vec{A}^T\vec{A}+\sum_{i=0}^p \vec{R}_i^T\vec{R}_i\right)\vec{x}^*
&= \lambda \vec{A}^T\vec{d}+\sum_{i=0}^p \vec{R}_i^T\vec{R}_i\bm{\mu}_i+
\vec{R}_\text{post}^T\bm{\xi}.
\end{align}

Inserting \eqref{eq:postprec} and solving for $\vec{x}^*$ we find
\begin{align}
\vec{x}^* &= \left(\vec{R}_\text{post}^T\vec{R}_\text{post}\right)^{-1} \left(  \lambda \vec{A}^T\vec{d} + \sum_{i=0}^p \vec{R}_i^T\vec{R}_i \bm{\mu}_i \right) +\left(\vec{R}_\text{post}^T\vec{R}_\text{post}\right)^{-1} 
\vec{R}_\text{post}^T
\bm{\xi} \nonumber\\
    &= \bm{\mu}_\text{post} + \vec{R}_\text{post}^{\dagger} \bm{\xi}, 
\end{align}
where $\vec{R}_\text{post}^{\dagger} = (\vec{R}_\text{post}^T\vec{R}_\text{post})^{-1}\vec{R}_\text{post}^T$ is the Moore--Penrose pseudoinverse. Note that $\vec{x}^*$ is a linear transformation of the Gaussian random variable $\bm{\xi}$. As $\vec{R}_0$ has full column rank, $\vec{R}_\text{post}^{\dagger}$ also has full row rank. It follows that $\vec{x}^*$ is Gaussian distributed with mean and covariance \cite{Bickel2001}:
\begin{subequations}
\begin{align}
    \mathrm{E}(\vec{x}^*) &= \bm{\mu}_\text{post} +  \vec{R}_\text{post}^{\dagger} \vec{0} = \bm{\mu}_\text{post}\\
    \mathrm{Cov}(\vec{x}^*) &= \vec{R}_\text{post}^{\dagger} \vec{I} (\vec{R}_\text{post}^{\dagger})^T 
    = \left(\vec{R}_\text{post}^T\vec{R}_\text{post}\right)^{-1}.
\end{align}
\end{subequations}

This confirms that $\vec{x}^*$ is distributed according to the target posterior distribution. That is, by solving \eqref{eq:post_linear_system} we obtain samples of CT reconstruction images distributed according to \eqref{eq:GaussianPosteriorDistribution}.

\section{Numerical Experiments}\label{sec:experiments}

FORCE Technology\footnote{\href{https://forcetechnology.com/}{forcetechnology.com}} has provided CT data from a deep-sea oil pipe acquired in a laboratory experiment using the offset fan-beam scan geometry illustrated in Figure \ref{fig:scan_geom}. We conduct synthetic and real data CT experiments to assess the performance of the proposed SGPs for obtaining useful CT reconstructions and relevant uncertainty estimates. The results are compared to a conventional deterministic reconstruction method where \eqref{eq:forward} is solved with the conjugate gradient least-squares (CGLS) algorithm \cite[Ch. 7]{CGLS} and stopped at semi-convergence. We interpret the posterior mean as the CT reconstruction and the posterior $95\%$ credible intervals for each pixel as associated uncertainty estimates. 

\subsection{Data description}

\subsubsection{Phantom}

\begin{figure}
    \centering
    \includegraphics[width=0.4\textwidth, trim={2.5cm 0.9cm 0cm 1.1cm},clip]{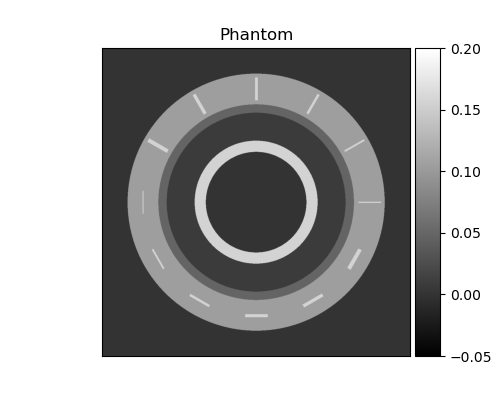}
    \includegraphics[width = 0.55\textwidth]{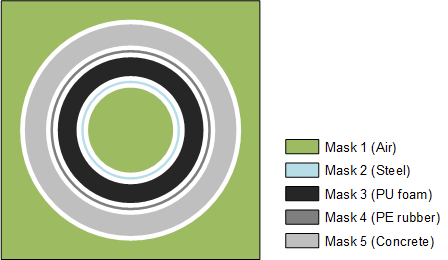}
    \caption{\textit{Left:} Pipe phantom modelling the real pipe seen in Figure \ref{fig:scan_geom}. \textit{Right:} SGP masks for $i=1,\ldots,5$ that each represent a material with attenuation coefficient $\alpha_i$ given in Table \ref{tab:abscoeff}. The white regions do not have any prior promoting specific attenuation coefficients. The index $i=0$ is used to denote the prior on the whole image domain (background and pipe).}
    \label{fig:SGPmask_synth}
\end{figure}

\begin{table}
    \caption{Best estimates for the material constants at beam intensity 2 MeV. }
    \label{tab:abscoeff}
    \begin{center}
    \begin{tabular}{lllllll}
    \toprule
    $i$ & Material & $\kappa [\text{cm}^2/\text{g}]$ & $\rho [\text{g}/\text{cm}^3]$ & $B$ & $w [\text{cm}]$& $\alpha [\text{cm}^{-1}]$ \\
    \midrule
    1 & Air & 0.044 & 0.0012 & - & - & $\sim 0$ \\
    2 & Steel\textsuperscript{a} & 0.042 & 7.9 & 2.013 & 4 & 0.16 \\
    3 & PU foam & 0.051\textsuperscript{b} & 0.15 & - & - & 0.0077 \\
    4 & PE rubber & 0.051 & 0.94 & - & - & 0.048 \\
    5 & Concrete & 0.046 & 2.3 & - & - & 0.11 \\
    \bottomrule
    \end{tabular}
    \end{center}
    \tabnote{\textsuperscript{a}Subject to build-up. \textsuperscript{b}Estimated by PE rubber's $\kappa$ due to missing value in database.}
\end{table}

Based on the known pipe structure and materials given in Figure \ref{fig:SGPmask_real}, we construct a phantom and corresponding SGP mask to perform synthetic experiments (see Figure \ref{fig:SGPmask_synth}). To avoid inverse crime, the phantom is defined on a $1024\times1024$ pixel image and the reconstruction is performed on a $512\times512$ pixel image that represents a $55 \times 55$ cm domain leading to pixels of size $\sim 1\times 1$ mm. The linear attenuation coefficients in each pipe layer are computed as described in Section \ref{sec:abscoeff}, with material constants from the NIST database \cite{NIST} and buildup factor from \cite{buildup}, assuming a monochromatic beam of energy 2 MeV\footnote{The beam energy is in fact a spectrum with 2 MeV mean energy.}. This beam intensity is an estimate, but around these energies the mass attenuation coefficients do not vary significantly. Therefore, we expect the computed linear attenuation coefficients are good estimates and they are listed in Table \ref{tab:abscoeff}. 

Note, that we add steel inclusions in the phantom's concrete layer to simulate steel reinforcement bars in the real pipe. The inclusions are in the radial or tangential directions and their widths are varied from 2--7 mm. These inclusions will not be accounted for in the SGP prior; this in order to test the prior's ability to reconstruct unexpected features. Furthermore, with the offset scan geometry, no X-rays will be parallel to near-radial edges, and therefore we expect that the radial inclusions are more difficult to reconstruct than the tangential (see, e.g., \cite{Riis2018} for more details).

\begin{figure}
    \centering
    \includegraphics[width = 0.49\textwidth]{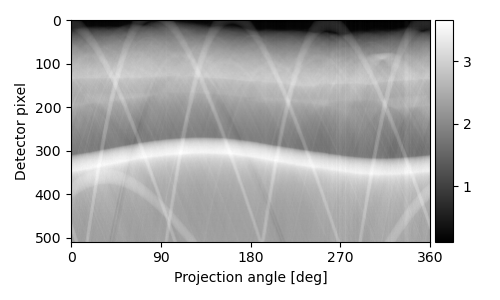} 
    \includegraphics[width = 0.49\textwidth]{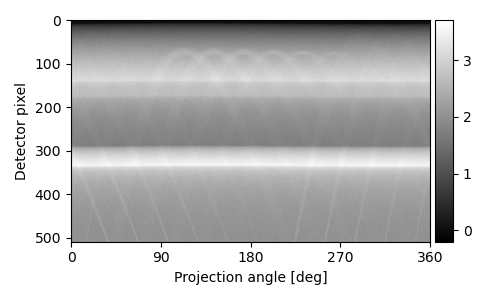}
    \vspace{-0.5cm}
    \caption{\textit{Left:} Real subsea pipe sinogram. \textit{Right:} Synthetic sinogram  with 2\% added noise.}
    \label{fig:sino}
\end{figure}
\subsubsection{Sinograms}

Using the pipe phantom, offset acquisition geometry, and 2\% noise level (i.e. $\lVert\vec{e}\rVert_2/\lVert\vec{Ax}\rVert_2 \sim 2\%$), we simulate a sinogram with 360 view angles distributed equally over a full rotation. Assuming the real sinogram has an equivalent noise distribution, we compare the synthetic and real data (Figure \ref{fig:sino}). Five bands in the sinograms mark the centre region, steel, PE, PU, and concrete layers in that order from the bottom of the sinograms. The bands in the real data sinogram follow sinusoidal curves because the pipe center did not coincide with the rotation center, and the phases of the sinusoidal bands vary because the pipe layers are not completely concentric. The overall similarity between the real data and synthetic data sinograms indicates that our phantom is realistic and the scan geometry close to the real setup. For sparse view angle experiments, we simulate data with only 50\%, 20\% or 10\% of the total view angles in the sinograms. Note that, since we reconstruct on a $512\times 512$ pixel image, the CT problem is under-determined in all our experiments.

\subsection{Numerical implementation}

\subsubsection{Bayesian parameter tuning} \label{sec:param_tune}

The posterior in \eqref{eq:GaussianPosterior0} requires choosing several parameters. The data noise of 2\% corresponds to a precision parameter $\lambda \approx 400$. The SGP includes several IID Gaussians for which we must determine mean and precision parameters. The means are defined as the estimated linear attenuation coefficients $\alpha_i$ for each pipe material reported in Table \ref{tab:abscoeff}. The IID Gaussian precision parameters are chosen heuristically such that the priors cover reasonable intervals around the estimated linear attenuation coefficients, representing some uncertainty in the estimates. We expect steel reinforcement bars in the concrete layer, which are not modeled in the prior, so the precision parameter is reduced here to allow a large range of possible values. Hence, the chosen precision parameters are: $\delta_1 = \delta_2= \delta_3 = \delta_4 = 1000$ and $\delta_5 = 500$. The final parameter choice is the GMRF precision $\delta_0$. If we fix all other precision parameters, we can perform a parameter sweep for $\delta_0$ values and choose the value that minimizes the root-mean-square-error (RMSE) with respect to the synthetic ground truth. An example of $\delta_0$ tuning can be seen in Figure \ref{fig:param_sweep} (left). We transfer the GMRF precision choices from the synthetic to the real data case, and check qualitatively that the corresponding reconstructions are acceptable.

\subsubsection{Sampling algorithm}\label{sec:algorithm}

We use the CGLS algorithm to solve \eqref{eq:postmean0} for MAP estimates of the posterior, and to solve \eqref{eq:post_linear_system} for realizations of the posterior. With this approach we avoid forming and solving the normal equations. Furthermore, this algorithm supports matrix-free implementation of $\vec{A}$ as required. We use the ASTRA toolbox to define the CT system and perform matrix-free, efficient, and GPU accelerated forward- and backprojections \cite{vanAarle2016,vanAarle2015,Palenstijn2011}. 

If the CGLS algorithm is run until convergence when solving \eqref{eq:post_linear_system}, we obtain independent realizations of the posterior \eqref{eq:GaussianPosteriorDistribution}. However, in order to reduce computation time, we apply a warm-start that uses the previous posterior realization as a starting guess for the next one. We fix the number of iterations to 10 per realization. This approach makes a burn-in phase necessary, i.e. some initial samples are discarded. Furthermore, the independence of the samples might be compromised, since we do not check for convergence. Therefore, we monitor the integrated autocorrelation time (IACT) of the samples \cite[p. 99]{Bardsley2018} to check for potential correlations. The idea is that, if $\text{IACT}\approx 1$, the computed samples are approximately independent.

\subsection{Synthetic experiments}

\subsubsection{Example of Bayesian inversion with SGP prior}

\begin{figure}
\centering
  \includegraphics[width=0.48\textwidth, trim={0cm 0.4cm 0cm 0cm},clip]{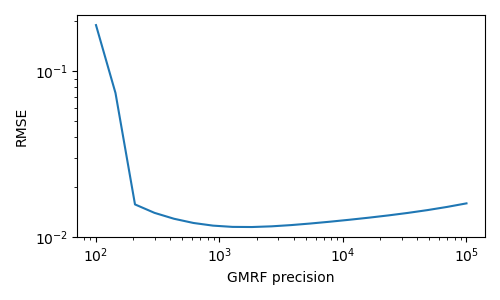}
  \includegraphics[width=0.48\textwidth, trim={0cm 0.4cm 0cm 0cm},clip]{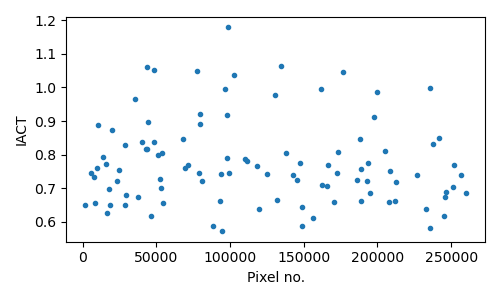}
  \caption{\textit{Left:} Example of GMRF precision parameter sweep from the synthetic experiment. From this, we choose $\delta_0=1000$ near the minimum. \textit{Right:} Examples of integrated autocorrelation time (IACT) for 100 randomly selected pixel chains.\label{fig:param_sweep}}
\end{figure}

We consider a sparse-angle case with only 72 equi-angular projections out of the original 360 (i.e. 20\%). We employ the SGP-F prior in the Bayesian inverse problem. First, we tune the GMRF precision parameter as described above in Section \ref{sec:param_tune}. Figure \ref{fig:param_sweep} (left) illustrates the parameter sweep from which we choose $\delta_0 = 1000$ since the RMSE is near a minimum. Then we compute 3000 samples from the posterior as described in Section \ref{sec:algorithm} and discard the first 1000 as burn-in. Since the reconstructed parameter $\vec{x}$ is high-dimensional, we compute the IACT at 100 randomly selected pixel chains and consider this as representative for the full parameter. Figure \ref{fig:param_sweep} (right) shows that IACT is close to or below one in all cases, indicating that the samples are quasi-independent.

\begin{figure}
    \centering
    \includegraphics[width = 0.345\textwidth, trim={4cm 4.5cm 11.5cm 0.2cm},clip]{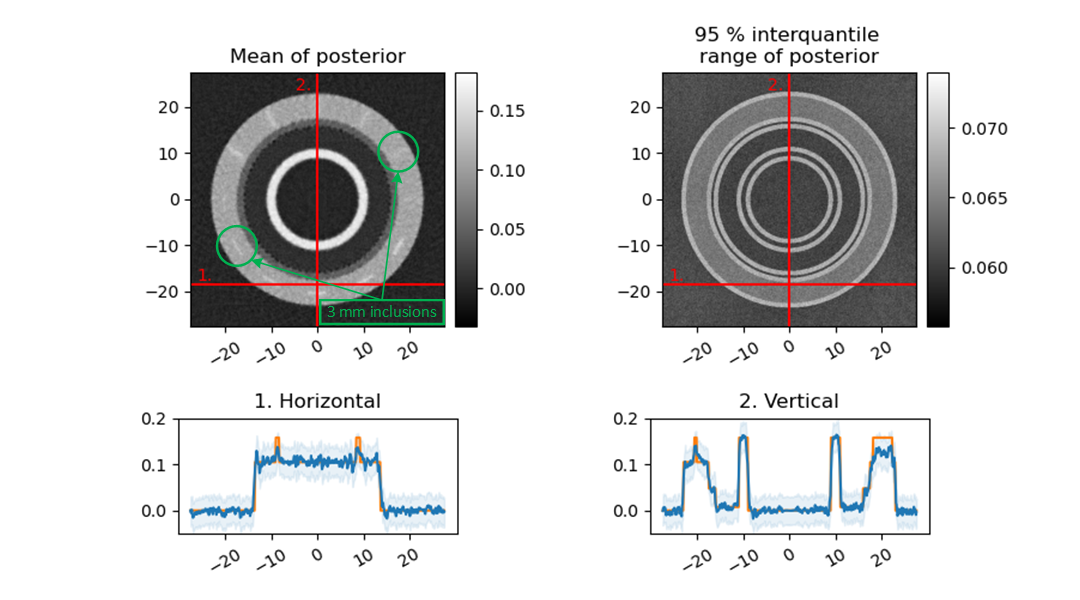}
    \includegraphics[width = 0.345\textwidth, trim={14cm 4.5cm 1.5cm 0.2cm},clip]{figures/20percent_SP01234/Slice_plot.png}
    \begin{minipage}{0.29\textwidth}
    \vspace{-5.9cm}
    \includegraphics[width =\textwidth, trim={2.8cm 0.5cm 13cm 8.2cm},clip]{figures/20percent_SP01234/Slice_plot.png}\\
    \includegraphics[width =\textwidth, trim={12.8cm 0.5cm 3cm 8.2cm},clip]{figures/20percent_SP01234/Slice_plot.png}
    \end{minipage}
    \vspace{-0.9cm}
    \caption{Synthetic test case with 20\% view angles and SGP-F prior. \textit{Left, center:} 2D images of the posterior mean and 95\% interquantile range. \textit{Right:} 1D slices of posterior mean (blue line)  and 95\% credibility interval (blue shade), compared to ground truth (orange).}
    \label{fig:1D_slice}
\end{figure}

Figure \ref{fig:1D_slice} shows the posterior mean, 95\% interquantile range (defined as the difference between the $2.5\%$ and $97.5\%$ quantiles), and 1D slices through the mean with 95\% credibility interval. We observe that the pipe layers are well reconstructed. However, the contrast of the steel inclusions in the concrete layer are dampened due to the choice of prior. The prior promotes the expected linear attenuation coefficient of the concrete, which is lower than that of steel. We also note higher noise in the concrete layer than in the remaining parts of the reconstruction due to $\delta_5<\delta_{1,\ldots,4}$. The 95\% interquantile range is a direct reflection of the SGP prior action. We see that imposing IID priors in the pipe layers lowers the uncertainty estimate of the posterior compared to areas without IID priors on them. Furthermore, priors with high precision values lead to lower uncertainty in the posterior. As expected, we also see that the radial steel inclusions in the concrete layer are more difficult to reconstruct than the tangential inclusions. For instance, the 3 mm radial inclusion is barely distinguishable, while the 3 mm tangential inclusion is fairly clear.

\setlength\tabcolsep{1.5pt}
\begin{figure}[ht!]
    \centering
    \begin{tabular}{c c c c c}
    & \textbf{Deterministic} & \textbf{GMRF} & \hspace{-4mm}\textbf{SGP-BG}\hspace{-4mm} & \textbf{SGP-F} \\
        \rotatebox[origin=l]{90}{\textbf{360 view angles}}
        & \includegraphics[width=0.23\textwidth, trim={2.5cm 1cm 2.2cm 1.1cm},clip]{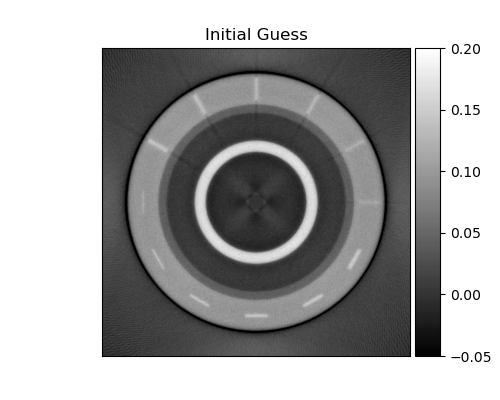}
        & \includegraphics[width=0.23\textwidth,trim={2.5cm 1cm 2.2cm 1.1cm},clip]{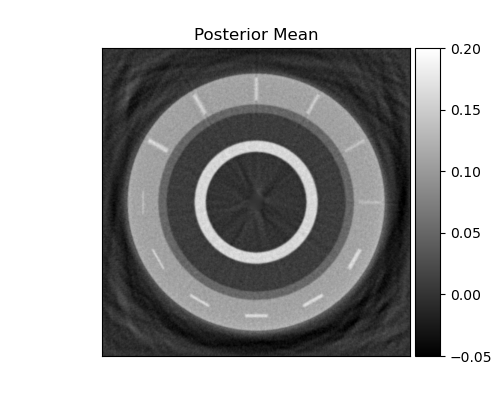} 
        & \includegraphics[width=0.23\textwidth,trim={2.5cm 1cm 2.2cm 1.1cm},clip]{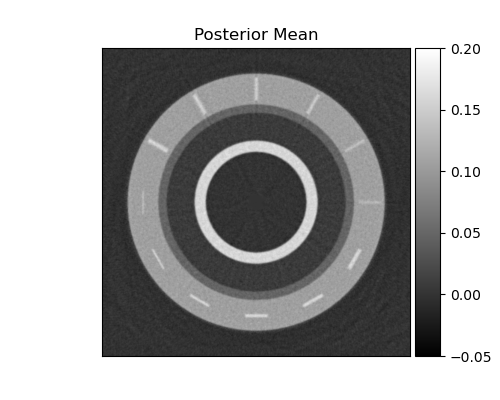}
        & \includegraphics[width=0.23\textwidth,trim={2.5cm 1cm 2.2cm 1.1cm},clip]{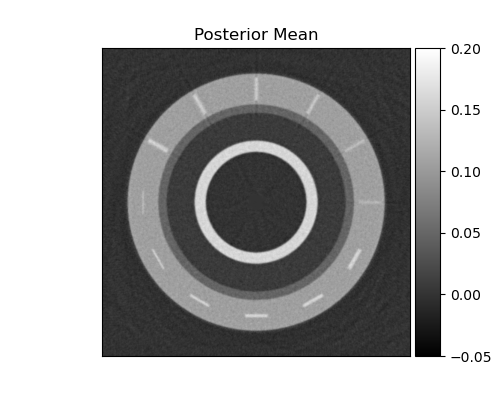} \\[-1mm]
        & \footnotesize (RMSE = 19.5$\times10^{-3}$) & \footnotesize (RMSE = 14.8$\times10^{-3}$) & \footnotesize (RMSE = 9.22$\times10^{-3}$) & \footnotesize (RMSE = 9.14$\times10^{-3}$) \\[1mm]
        \rotatebox[origin=l]{90}{\textbf{50\% view angles}}
        & \includegraphics[width=0.23\textwidth, trim={2.5cm 1cm 2.2cm 1.1cm},clip]{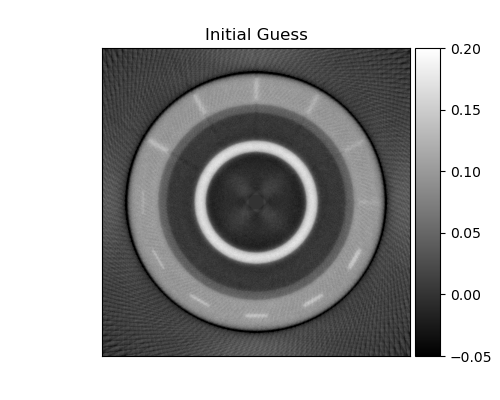}
        & \includegraphics[width=0.23\textwidth,trim={2.5cm 1cm 2.2cm 1.1cm},clip]{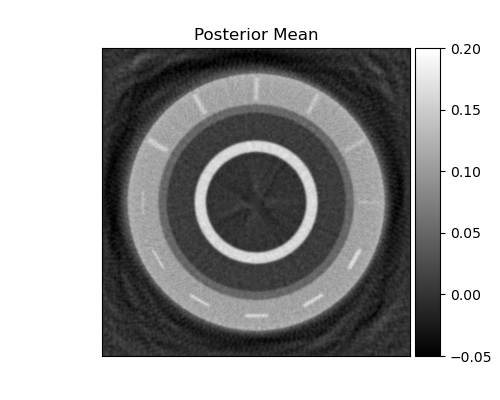} 
        & \includegraphics[width=0.23\textwidth,trim={2.5cm 1cm 2.2cm 1.1cm},clip]{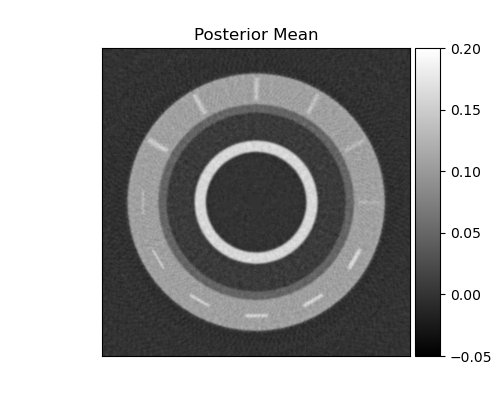} 
        & \includegraphics[width=0.23\textwidth,trim={2.5cm 1cm 2.2cm 1.1cm},clip]{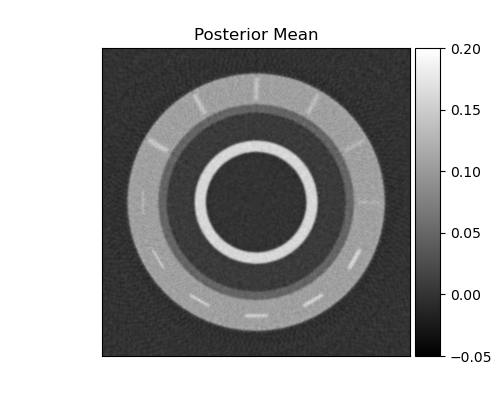}\\[-1mm]
        & \footnotesize(RMSE = 20.6$\times10^{-3}$) & \footnotesize(RMSE = 16.7$\times10^{-3}$) & \footnotesize(RMSE = 10.1$\times10^{-3}$) & \footnotesize(RMSE = 9.91$\times10^{-3}$) \\[1mm]
        \rotatebox[origin=l]{90}{\textbf{20\% view angles}}
        & \includegraphics[width=0.23\textwidth, trim={2.5cm 1cm 2.2cm 1.1cm},clip]{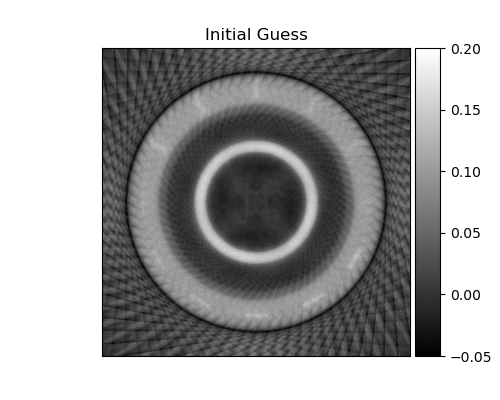}
        & \includegraphics[width=0.23\textwidth,trim={2.5cm 1cm 2.2cm 1.1cm},clip]{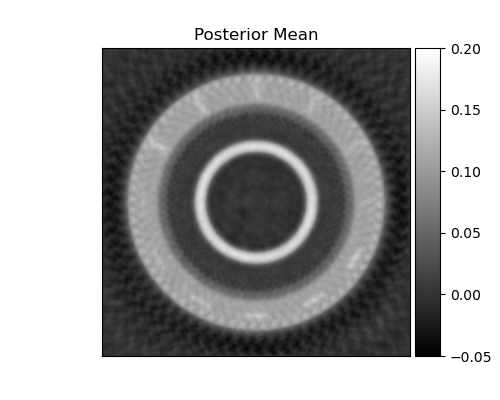} 
        & \includegraphics[width=0.23\textwidth,trim={2.5cm 1cm 2.2cm 1.1cm},clip]{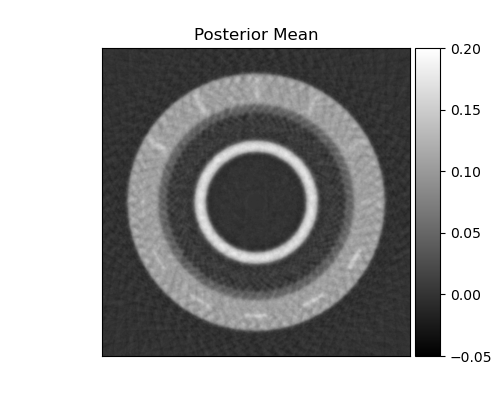} 
        & \includegraphics[width=0.23\textwidth,trim={2.5cm 1cm 2.2cm 1.1cm},clip]{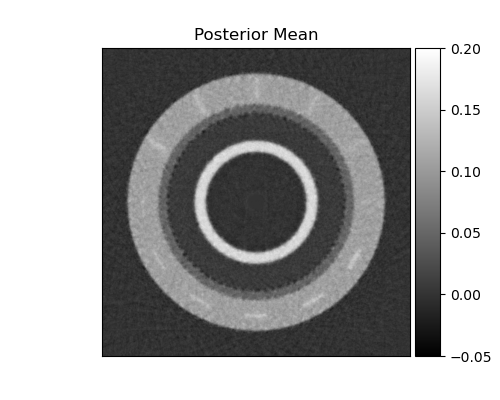} \\[-1mm]
        & \footnotesize(RMSE = 25.1$\times10^{-3}$) & \footnotesize(RMSE = 16.5$\times10^{-3}$) & \footnotesize(RMSE = 12.2$\times10^{-3}$) & \footnotesize(RMSE = 11.6$\times10^{-3}$) \\[1mm]
        \rotatebox[origin=l]{90}{\textbf{10\% view angles}}
        & \includegraphics[width=0.23\textwidth, trim={2.5cm 1cm 2.2cm 1.1cm},clip]{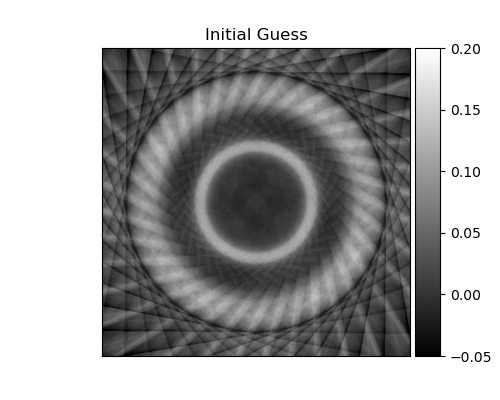}
        & \includegraphics[width=0.23\textwidth,trim={2.5cm 1cm 2.2cm 1.1cm},clip]{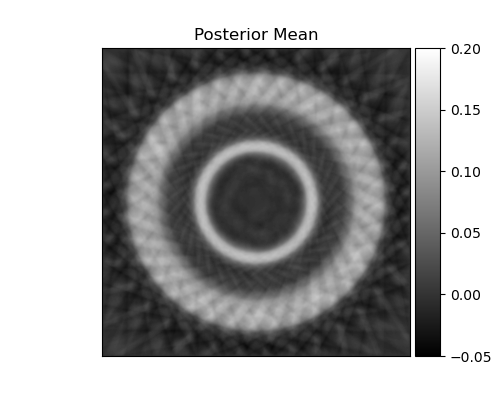} 
        & \includegraphics[width=0.23\textwidth,trim={2.5cm 1cm 2.2cm 1.1cm},clip]{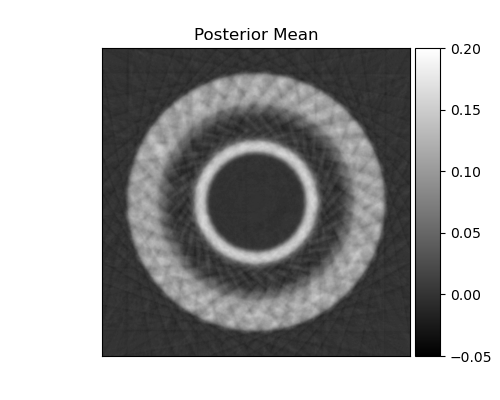} 
        & \includegraphics[width=0.23\textwidth,trim={2.5cm 1cm 2.2cm 1.1cm},clip]{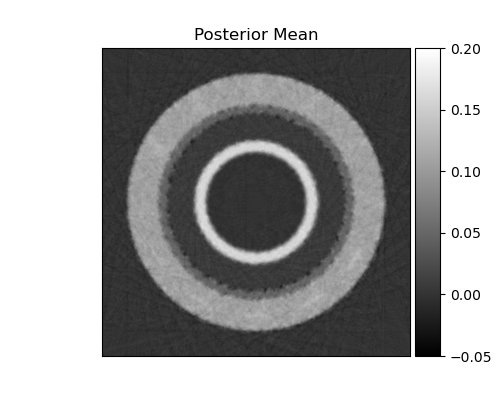} \\[-1mm]
        & \footnotesize(RMSE = 30.6$\times10^{-3}$) & \footnotesize(RMSE = 18.9$\times10^{-3}$) & \footnotesize(RMSE = 15.0$\times10^{-3}$) & \footnotesize(RMSE = 12.5$\times10^{-3}$)
    \end{tabular}
    \caption{Posterior means and RMSE with increasingly informative SGPs (in columns) and decreasing number of view angles (in rows). The color range is $[-0.05, 0.20]$ for all images.}
    \label{fig:synth_post_means}
\end{figure}

\setlength\tabcolsep{1.5pt}
\begin{figure}[ht!]
    \centering
    \begin{tabular}{c c c c c}
    & \textbf{GMRF} \hspace{6mm} & \hspace{-2mm}\textbf{SGP-BG}\hspace{2mm} & \hspace{-4mm}\textbf{SGP-F} \hspace{6mm}\\
        \rotatebox[origin=l]{90}{\textbf{360 view angles}}
        & \includegraphics[width=0.27\textwidth,trim={2.5cm 1cm 0cm 1.1cm},clip]{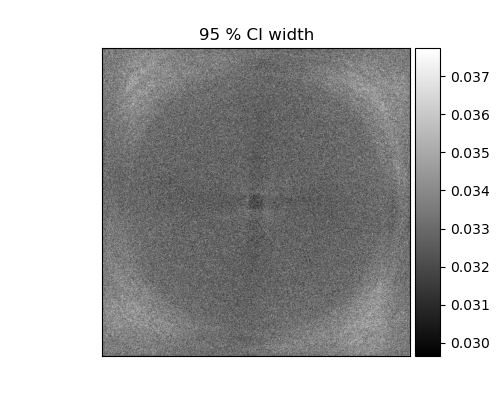} 
        & \includegraphics[width=0.27\textwidth,trim={2.5cm 1cm 0cm 1.1cm},clip]{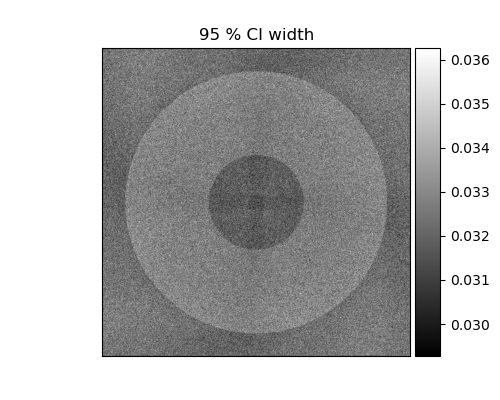} 
        & \includegraphics[width=0.27\textwidth,trim={2.5cm 1cm 0cm 1.1cm},clip]{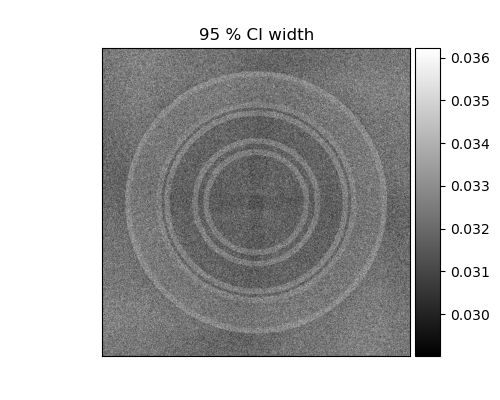}\\[-1mm]
        & \footnotesize($\delta_0 = 4000$) \hspace{6mm} & \footnotesize($\delta_0 = 4000$)\hspace{6mm} & \footnotesize($\delta_0 = 4000$)\hspace{6mm} \\
        \rotatebox[origin=l]{90}{\textbf{50\% view angles}}
        & \includegraphics[width=0.27\textwidth,trim={2.5cm 1cm 0cm 1.1cm},clip]{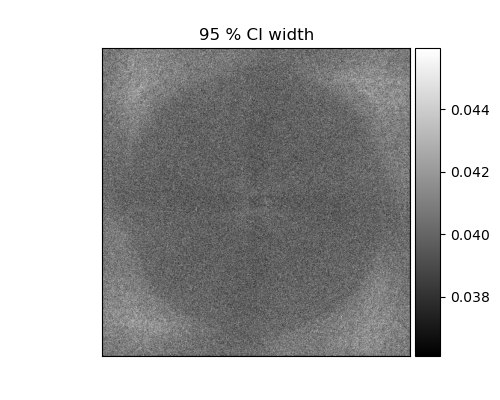} 
        & \includegraphics[width=0.27\textwidth,trim={2.5cm 1cm 0cm 1.1cm},clip]{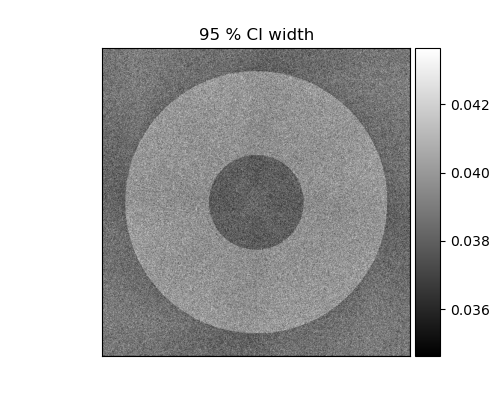} 
        & \includegraphics[width=0.27\textwidth,trim={2.5cm 1cm 0cm 1.1cm},clip]{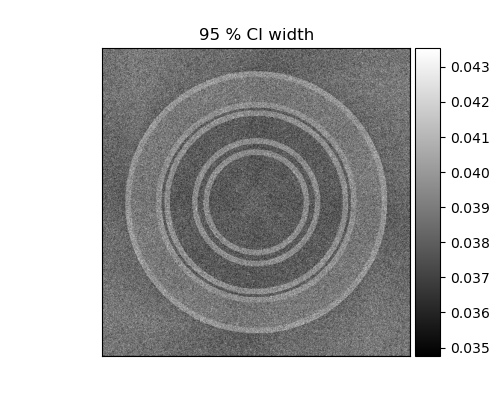}\\[-1mm]
        & \footnotesize($\delta_0 = 3000$)\hspace{6mm} & \footnotesize($\delta_0 = 3000$)\hspace{6mm} & \footnotesize($\delta_0 = 3000$)\hspace{6mm} \\
        \rotatebox[origin=l]{90}{\textbf{20\% view angles}}
        & \includegraphics[width=0.27\textwidth,trim={2.5cm 1cm 0cm 1.1cm},clip]{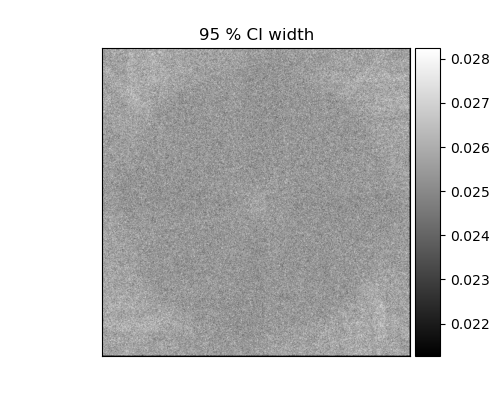} 
        & \includegraphics[width=0.27\textwidth,trim={2.5cm 1cm 0cm 1.1cm},clip]{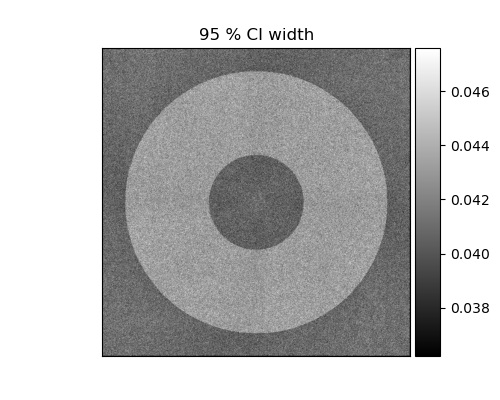} 
        & \includegraphics[width=0.27\textwidth,trim={2.5cm 1cm 0cm 1.1cm},clip]{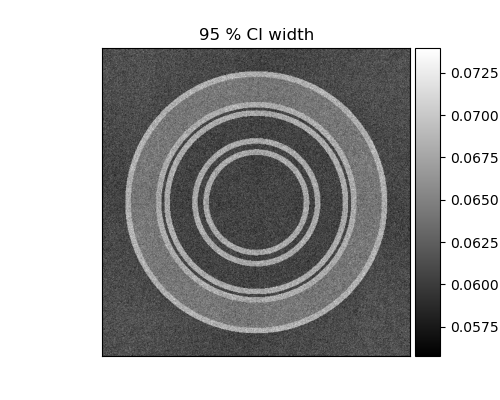} \\[-1mm]
        & \footnotesize($\delta_0 = 10000$)\hspace{6mm} & \footnotesize($\delta_0 = 3000$)\hspace{6mm} & \footnotesize($\delta_0 = 1000$)\hspace{6mm} \\
        \rotatebox[origin=l]{90}{\textbf{10\% view angles}}
        & \includegraphics[width=0.27\textwidth,trim={2.5cm 1cm 0cm 1.1cm},clip]{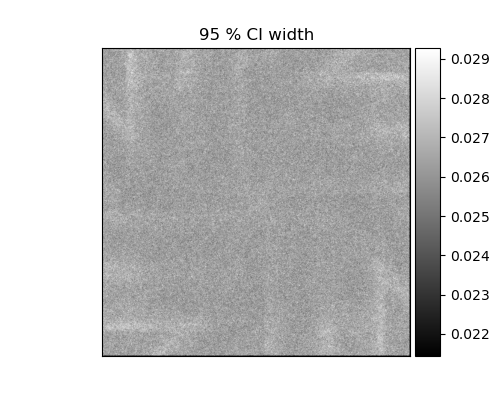} 
        & \includegraphics[width=0.27\textwidth,trim={2.5cm 1cm 0cm 1.1cm},clip]{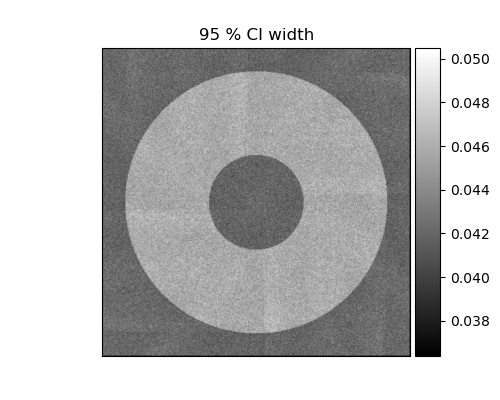} 
        & \includegraphics[width=0.27\textwidth,trim={2.5cm 1cm 0cm 1.1cm},clip]{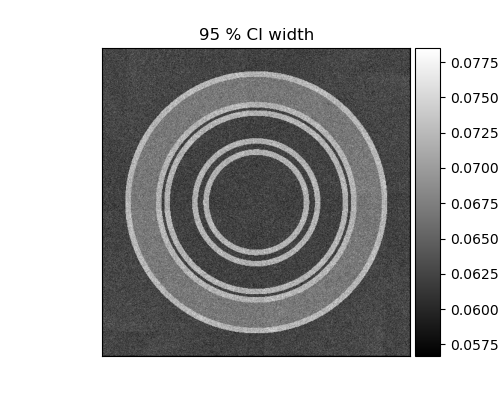}\\[-1mm]
        & \footnotesize($\delta_0 = 10000$)\hspace{6mm} & \footnotesize($\delta_0 = 3000$)\hspace{6mm} & \footnotesize($\delta_0 = 1000$) \hspace{6mm}
    \end{tabular}
    
    \caption{Posterior 95\% interquantile range with increasingly informative SGPs (in columns) and decreasing number of view angles (in rows).}
    \label{fig:synth_post_interq}
\end{figure}

\subsubsection{Comparative study}

We investigate the effect of the SGP when reducing the number of view angles. Figure \ref{fig:synth_post_means} shows posterior means on a grid, where the number of view angles decreases by rows, and the prior information increases by columns. We observe that, as we add more structural information, the reconstructions improve in terms of 1) reduced streak artifacts, 2) increased layer contrast, 3) reduced RMSE, and 4) detection of the steel inclusions that are hidden in artifacts. The action of the prior is especially clear for 10\% and 20\% view angles, indicating that the structural priors become more effective in sparse view angle cases.

The deterministic reconstruction with 360 view angles produces an image where all the main features of the pipe are visible. However, we also see artifacts: 1) along the radial steel inclusions, because these edges have no tangential rays, 2) in the centre of the image, because there is no data, and 3) a dark ring just outside of the pipe, which is a typical artifact that appears at the fan-beam edge. As we decrease the number of view angles, we get more dominating artifacts from the fan-beam edges. When we include the GMRF prior, the artifacts are smoothed out, but they are still very dominating in the image. Thus, this standard prior is not sufficient to produce good reconstructions, especially for sparse view angle settings. Alternatively, when using the SGP-BG prior, we manage to reduce the artifacts in the background substantially. Furthermore, we see increased contrast between pipe materials for the 10\% and 20\% view angle cases. This is due to the reduced precision parameter of the smoothing GMRF prior, made possible by the added information from the SGP-BG prior. When we use the SGP-F prior, we can suppress artifacts in the background and pipe layers. In the concrete layer, however, we note that unexpected inclusions (here steel reinforcements) that are not modeled in the prior, are also removed. Recall, that we did not impose any IID prior in boundary areas, and therefore the artifacts are still present there. This we see as a `ragged' boundary, especially for sparse view angle cases.

We see that, with the SGP-F prior, we can resolve tangential steel inclusions down to widths of 2 mm using only 20\% of the view angles. The radial inclusions are well reconstructed down to 4 mm. The smaller 2-3 mm radial inclusions are still visible but faint and blurred. Compared to the deterministic reconstruction with no prior, or to the reconstruction with a GMRF prior alone, the results still show an improvement. Note that, in those cases, the narrow radial steel inclusions disappear among the artifacts. Moreover, for the 10\% view angles case, it becomes difficult to reconstruct any of the steel inclusions using any of the priors. Only the widest tangential inclusions appear faintly in the reconstructions.

In Figure \ref{fig:synth_post_interq}, we show the 95\% interquantile ranges of the posteriors which reflect the uncertainty in the reconstructions. Note that the colorbar for each image is different. The GMRF prior yields a posterior uncertainty pattern that reflects the acquisition geometry and that comes from the uncertainty of the likelihood. For the SGP priors, the posterior uncertainty pattern reflects the structure of the prior. For SGP-F, we observe lowest uncertainty in the background and in the steel, PU, and PE layers (regions 1-4), where we find the highest prior precision parameters. The concrete layer (region 5) has a lower prior precision, and thus it reflects higher uncertainty. The highest uncertainty is found at the boundary regions, where no particular attenuation coefficient is promoted. We note that the boundary regions seem to have relatively higher uncertainty for 10\% and 20\% view angles, than for 50\% and 100\%. This is caused by the lower GMRF precision $\delta_0$ used in the two former cases. Finally, the likelihood structure is visible for the SGP-F prior cases with 50\% and 100\% view angles, but not with 10\% and 20\%. We suspect this is due to the prior having a relatively larger influence than the likelihood when data is sparse.

\setlength\tabcolsep{1.5pt}
\begin{figure}[ht!]
    \centering
    \begin{tabular}{c c c c c}
    & \textbf{Deterministic} & \textbf{GMRF} & \hspace{-4mm}\textbf{SGP-BG}\hspace{-4mm} & \textbf{SGP-F}\\
        \rotatebox[origin=l]{90}{\textbf{360 view angles}}
        & \includegraphics[width=0.23\textwidth, trim={2.5cm 1cm 2.2cm 1.1cm},clip]{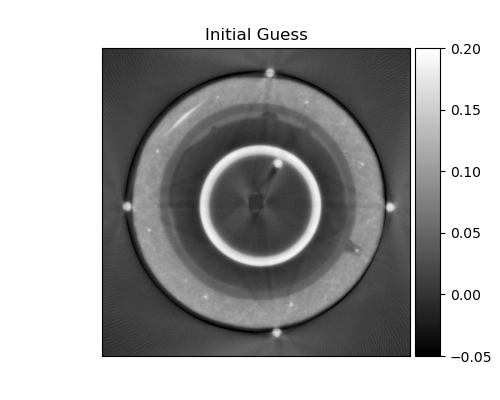}
        & \includegraphics[width=0.23\textwidth,trim={2.5cm 1cm 2.2cm 1.1cm},clip]{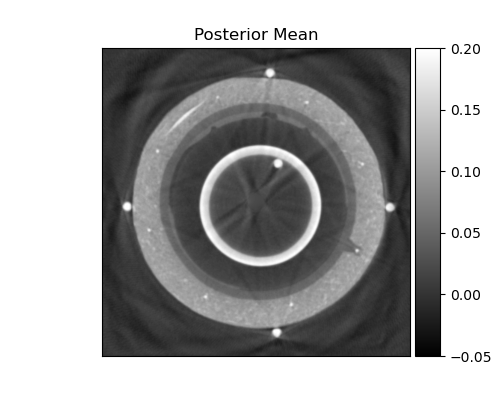} 
        & \includegraphics[width=0.23\textwidth,trim={2.5cm 1cm 2.2cm 1.1cm},clip]{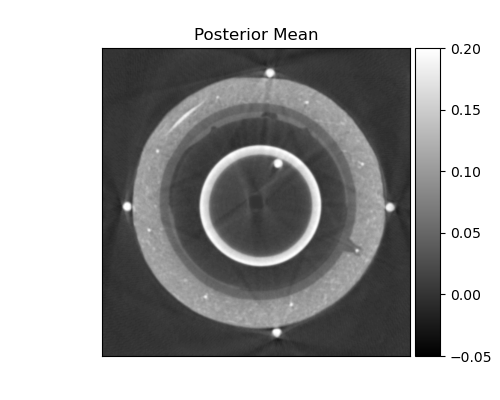} 
        & \includegraphics[width=0.23\textwidth,trim={2.5cm 1cm 2.2cm 1.1cm},clip]{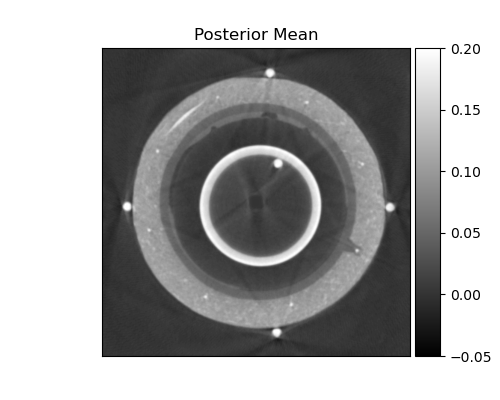} \\
        \rotatebox[origin=l]{90}{\textbf{50\% view angles}}
        & \includegraphics[width=0.23\textwidth, trim={2.5cm 1cm 2.2cm 1.1cm},clip]{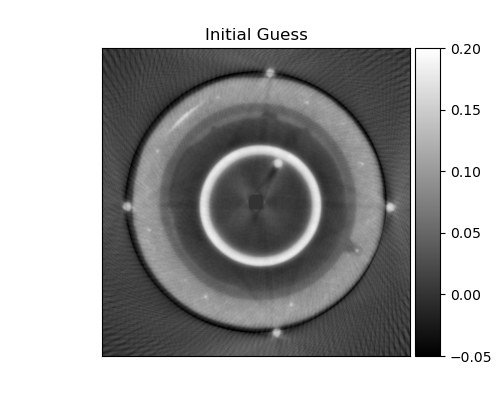}
        & \includegraphics[width=0.23\textwidth,trim={2.5cm 1cm 2.2cm 1.1cm},clip]{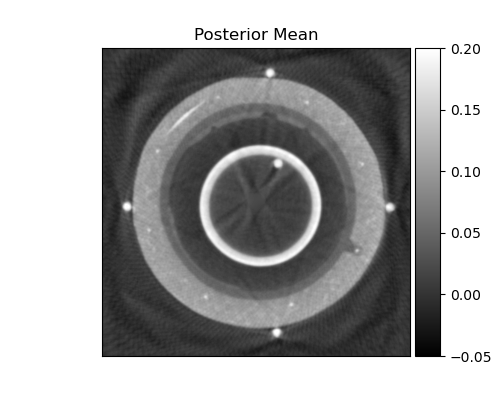} 
        & \includegraphics[width=0.23\textwidth,trim={2.5cm 1cm 2.2cm 1.1cm},clip]{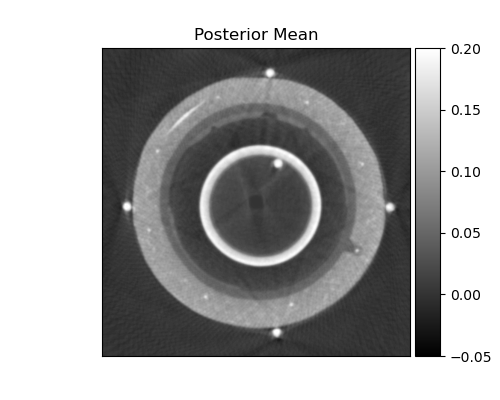} 
        & \includegraphics[width=0.23\textwidth,trim={2.5cm 1cm 2.2cm 1.1cm},clip]{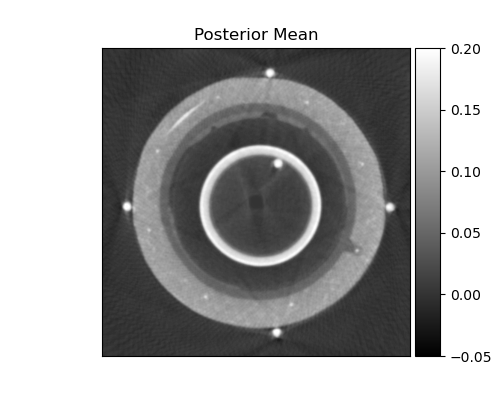} \\
        \rotatebox[origin=l]{90}{\textbf{20\% view angles}}
        & \includegraphics[width=0.23\textwidth, trim={2.5cm 1cm 2.2cm 1.1cm},clip]{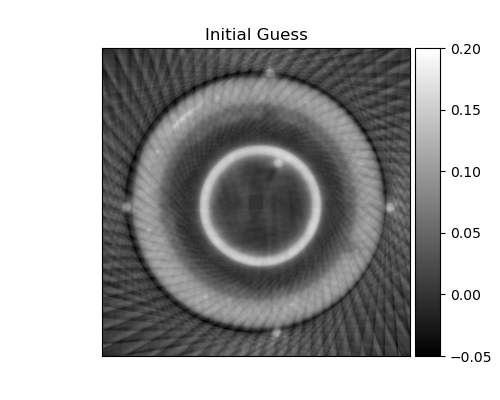}
        & \includegraphics[width=0.23\textwidth,trim={2.5cm 1cm 2.2cm 1.1cm},clip]{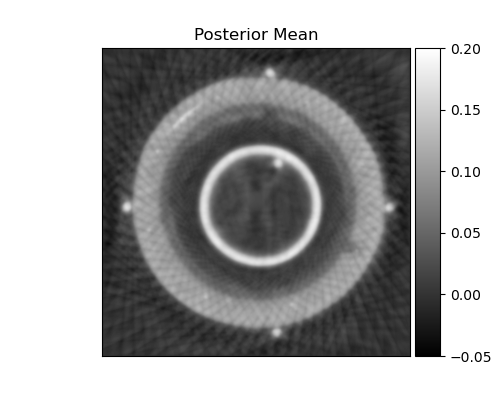} 
        & \includegraphics[width=0.23\textwidth,trim={2.5cm 1cm 2.2cm 1.1cm},clip]{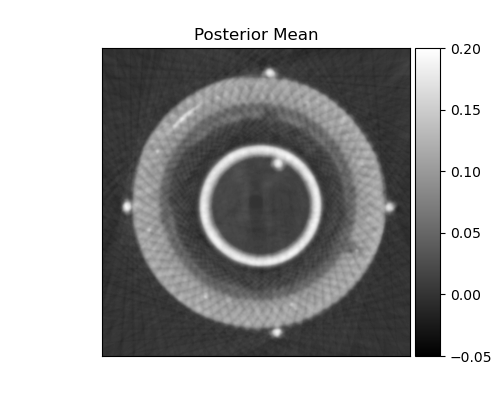} 
        & \includegraphics[width=0.23\textwidth,trim={2.5cm 1cm 2.2cm 1.1cm},clip]{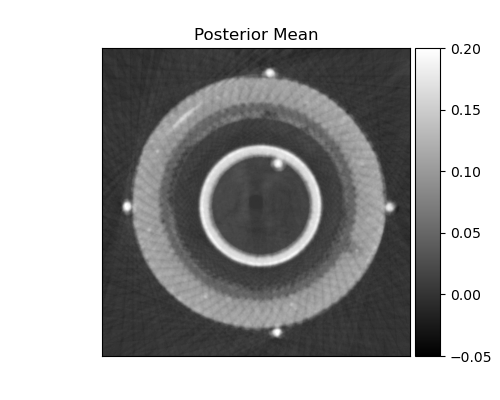} \\
        \rotatebox[origin=l]{90}{\textbf{10\% view angles}}
        & \includegraphics[width=0.23\textwidth, trim={2.5cm 1cm 2.2cm 1.1cm},clip]{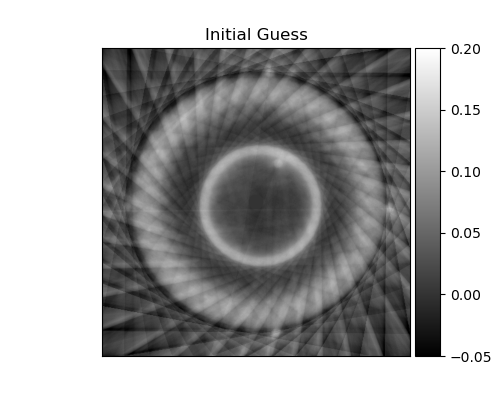}
        & \includegraphics[width=0.23\textwidth,trim={2.5cm 1cm 2.2cm 1.1cm},clip]{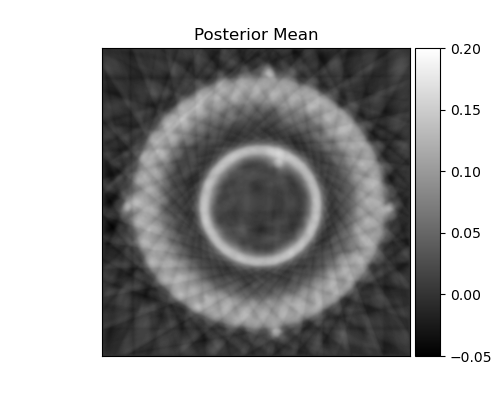} 
        & \includegraphics[width=0.23\textwidth,trim={2.5cm 1cm 2.2cm 1.1cm},clip]{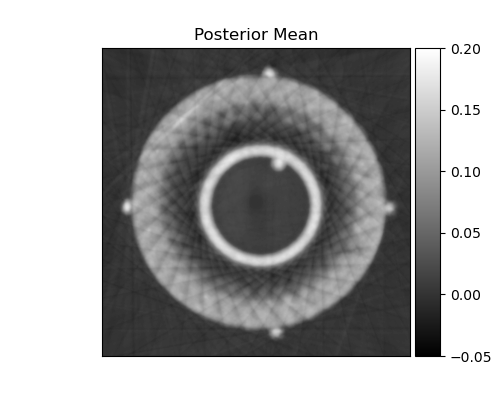} 
        & \includegraphics[width=0.23\textwidth,trim={2.5cm 1cm 2.2cm 1.1cm},clip]{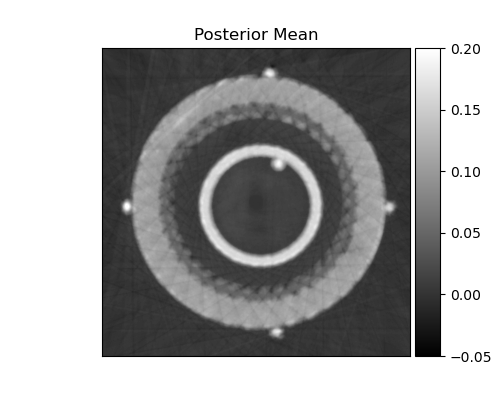} 
    \end{tabular}
    \caption{Posterior means with increasingly informative SGPs (in columns) and decreasing number of view angles (in rows). The color range is $[-0.05, 0.20]$ for all images.}
    \label{fig:real_post_means}
\end{figure}

\setlength\tabcolsep{1.5pt}
\begin{figure}[ht!]
    \centering
    \begin{tabular}{c c c c c}
    & \textbf{GMRF} \hspace{6mm} & \hspace{-2mm}\textbf{SGP-BG}\hspace{2mm} & \hspace{-4mm}\textbf{SGP-F} \hspace{6mm}\\
        \rotatebox[origin=l]{90}{\textbf{360 view angles}}
        & \includegraphics[width=0.27\textwidth,trim={2.5cm 1cm 0cm 1.1cm},clip]{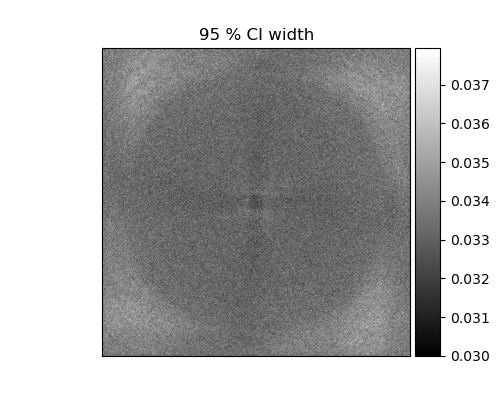} 
        & \includegraphics[width=0.27\textwidth,trim={2.5cm 1cm 0cm 1.1cm},clip]{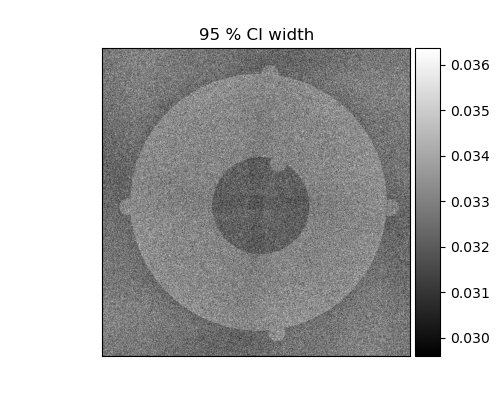} 
        & \includegraphics[width=0.27\textwidth,trim={2.5cm 1cm 0cm 1.1cm},clip]{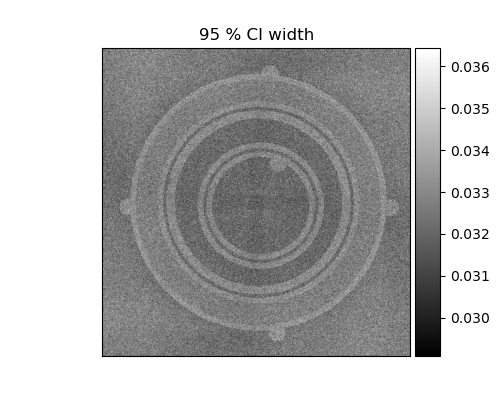}\\[-1mm]
        & \footnotesize($\delta_0 = 4000$) \hspace{6mm} & \footnotesize($\delta_0 = 4000$)\hspace{6mm} & \footnotesize($\delta_0 = 4000$)\hspace{6mm} \\
        \rotatebox[origin=l]{90}{\textbf{50\% view angles}}
        & \includegraphics[width=0.27\textwidth,trim={2.5cm 1cm 0cm 1.1cm},clip]{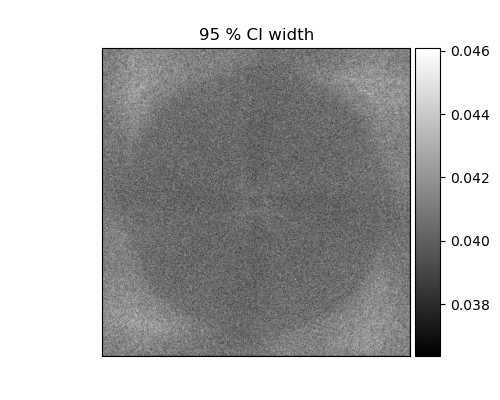} 
        & \includegraphics[width=0.27\textwidth,trim={2.5cm 1cm 0cm 1.1cm},clip]{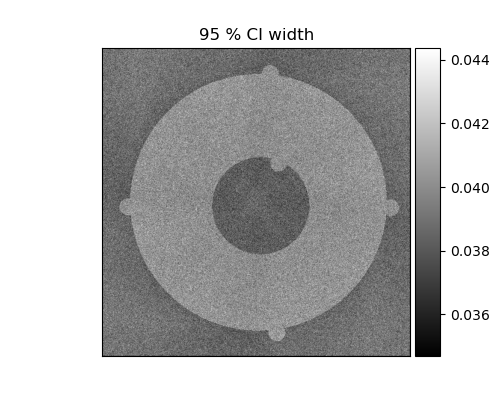} 
        & \includegraphics[width=0.27\textwidth,trim={2.5cm 1cm 0cm 1.1cm},clip]{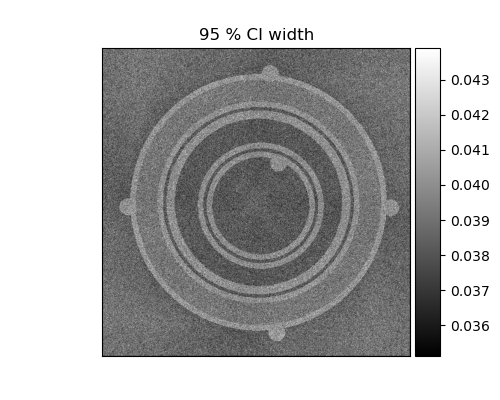}\\[-1mm]
        & \footnotesize($\delta_0 = 3000$)\hspace{6mm} & \footnotesize($\delta_0 = 3000$)\hspace{6mm} & \footnotesize($\delta_0 = 3000$)\hspace{6mm} \\
        \rotatebox[origin=l]{90}{\textbf{20\% view angles}}
        & \includegraphics[width=0.27\textwidth,trim={2.5cm 1cm 0cm 1.1cm},clip]{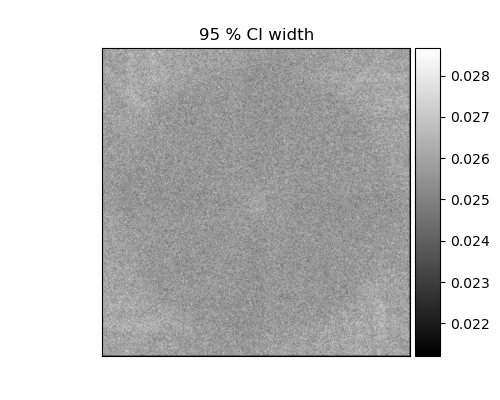} 
        & \includegraphics[width=0.27\textwidth,trim={2.5cm 1cm 0cm 1.1cm},clip]{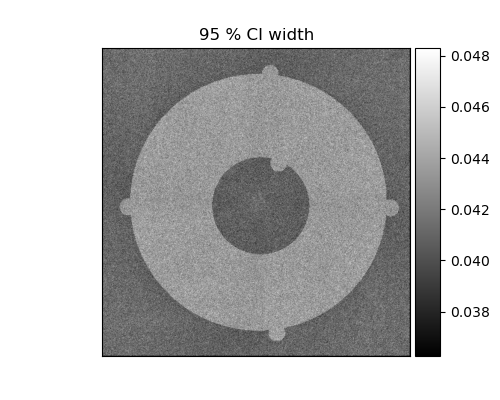} 
        & \includegraphics[width=0.27\textwidth,trim={2.5cm 1cm 0cm 1.1cm},clip]{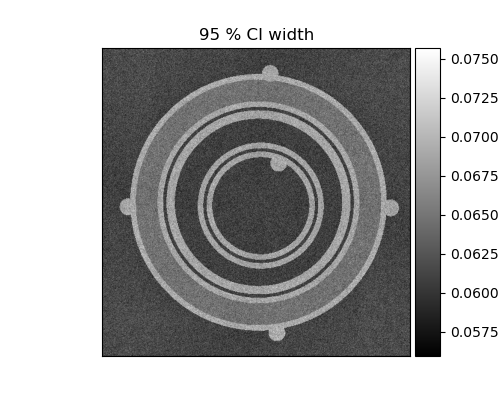} \\[-1mm]
        & \footnotesize($\delta_0 = 10000$)\hspace{6mm} & \footnotesize($\delta_0 = 3000$)\hspace{6mm} & \footnotesize($\delta_0 = 1000$)\hspace{6mm} \\
        \rotatebox[origin=l]{90}{\textbf{10\% view angles}}
        & \includegraphics[width=0.27\textwidth,trim={2.5cm 1cm 0cm 1.1cm},clip]{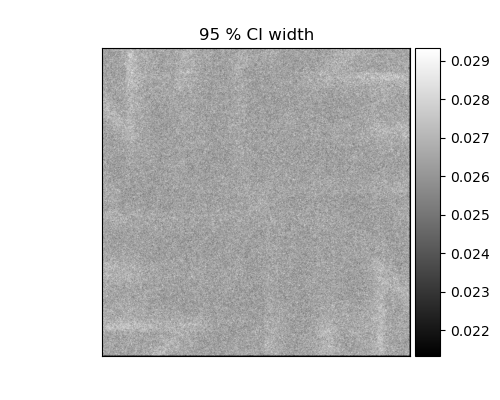} 
        & \includegraphics[width=0.27\textwidth,trim={2.5cm 1cm 0cm 1.1cm},clip]{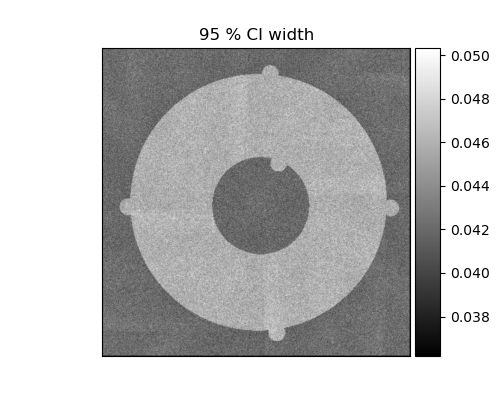} 
        & \includegraphics[width=0.27\textwidth,trim={2.5cm 1cm 0cm 1.1cm},clip]{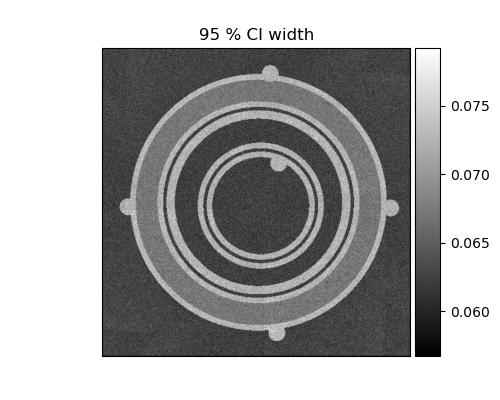}\\[-1mm]
        & \footnotesize($\delta_0 = 10000$)\hspace{6mm} & \footnotesize($\delta_0 = 3000$)\hspace{6mm} & \footnotesize($\delta_0 = 1000$) \hspace{6mm}
    \end{tabular}
    \caption{Real data: Posterior 95\% interquantile range with increasingly informative SGPs (in columns) and decreasing number of view angles (in rows).}
    \label{fig:real_post_interq}
\end{figure}

\subsection{Real data experiments}

We now apply the SGPs to real data and obtain reconstructions with image qualities that appear comparable to the synthetic experiments. As for the synthetic data, we will study the effect of the SGP priors in combination with sparse view angles data. Figure \ref{fig:real_post_means} shows the posterior means on a grid, where the priors and number of view angles are varied. Our observations in the synthetic data experiments are in good agreement with the real data case. Again, we see dominating artifacts when decreasing the number of view angles in the deterministic reconstruction with no prior information. We also see that they are reduced by introducing structural information. In the 20\% view angles case, we observe that the SGPs not only reduce artifacts, but also help increase contrast between materials. In the most extreme case, with only 10\% view angles, we lose a lot of details in the reconstruction, such as the irregularities between the PU and PE layers. We do, however, still see traces of the large steel inclusion in the concrete layer. 

In Figure \ref{fig:real_post_interq}, we see 95\% interquantile ranges associated with the posteriors. Again, we observe that the acquisition geometry is reflected when using the GMRF prior, and that the SGP structures are reflected (layers) when using these particular priors. Note that the steel rods used to fixate the pipe during measurements do not have any IID prior promoting a particular attenuation coefficient.

\section{Discussion and Conclusions} \label{sec:discussion}

We proposed a method of introducing prior structural information in Bayesian CT reconstruction. Structural Gaussian priors are constructed based on well-known geometry and materials to promote expected attenuation coefficients in CT reconstructions. We applied the SGP in a case study of CT imaging of a subsea pipe with an offset fan-beam acquisition geometry. Numerical experiments with real and synthetic data showed that SGPs improve reconstructions in terms of suppressed artifacts and enhanced contrast between materials, compared to traditional methods with no prior information or with a standard smoothing GMRF prior. The SGPs were especially powerful in setups with sparse view angles. Information about all layers of the pipe is required to construct the SGP-F. Using this prior yielded the best results in the numerical experiments. In some cases, we might only have access to the pipe's exterior dimensions. Then, the SGP-BG prior is applicable, and also in such cases, the numerical results showed substantial improvements of the reconstructions compared to using a GMRF prior or no prior. 

Since the posterior is Gaussian and defined in terms of square-root precision matrices, it can be efficiently sampled by solving a system of linear equations, yielding approximately independent samples. This is a useful property compared to conventional Markov chain Monte Carlo (MCMC) methods such as Metropolis Hastings \cite{Hastings1970} or more modern Hamiltonian MCMC methods \cite{Duane1987}, where independent samples cannot be guaranteed. To overcome this issue, MCMC methods often require a `thinning' process where many of the samples are discarded. Therefore, MCMC methods require many more simulations, and thus, much longer computation times to obtain a sufficiently large quasi-independent sample size.

If uncertainty estimates are not particularly required, our method based on SGP priors can also achieve reconstructions using MAP estimates that can be computed extremely efficiently. However, since one of our main goals is uncertainty quantification, we chose to sample the posterior. 

Interpretation of the posterior uncertainty estimates is difficult and requires more study. We observed that the GMRF prior gives rise to uncertainty estimates that reflect the acquisition geometry, and this is in agreement with the results in \cite{Bardsley2012b}. When imposing the SGP, we noticed, that the uncertainty estimates directly reflected the prior. The prior's effect on the uncertainty estimates has also been seen in other studies. For instance, edge-preserving priors give rise to high uncertainty on discontinuities in CT reconstructions \cite{Suuronen2020,Uribe2021}. In future work, we aim to further improve uncertainty estimates of the pipe and defect structure, for example to help distinguish artifacts from real features in the reconstruction. 

\section*{Acknowledgements}

We would like to thank FORCE Technology for providing real data and information regarding subsea pipes.

\section*{Funding}

This work was supported by The Villum Foundation (grant no. 25893).

\bibliographystyle{tfnlm}
\bibliography{iopart-num.bib}

\end{document}